\documentclass[a4paper]{amsart}
\usepackage{amssymb}

\input xy
\xyoption{all}

\newcommand{\pr}{\mathrm{pr}}

\newcommand{\id}{\mathrm{id}}

\newcommand{\R}{\mathbb R}

\newcommand{\N}{\mathbb N}

\newtheorem{theorem}{Theorem}
\newtheorem{df}{Definition}
\newtheorem{lemma}{Lemma}

\newtheorem{proposition}{Proposition}

\newtheorem{question}{Question}

\begin{document}

\title{On idempotent convexities and idempotent barycenter maps}


\author{Dawid Krasi\'{n}ski, Taras Radul}

\maketitle

Institute of Mathematics, Casimir the Great University of Bydgoszcz, Poland;
\newline
Department of Mechanics and Mathematics, Ivan Franko National University of Lviv,
Universytettska st., 1. 79000 Lviv, Ukraine.
\newline
e-mail: tarasradul@yahoo.co.uk



\begin{abstract} We consider an isomorphism between the  idempotent convexity based on the maximum and the addition operations and the idempotent measure convexity on the maximum and the multiplication  operations. We use this isomorphism to investigate topological properties of the barycenter map related to the maximum and the multiplication  operations.
 \end{abstract}

\maketitle

\section{Introduction} Idempotent mathematics has been a rapidly growing field in the last few decades. In traditional mathematics, addition is a fundamental operation. However, in idempotent mathematics, a different operation, typically the maximum operation, replaces addition.  Another classical multiplication operation was changed by different operations: addition, minimum, t-norm etc, depending on the specific model needed for a given application. These changes have led to numerous applications in various fields, including mathematics, mathematical physics, game theory, computer science, and economics. The versatility of idempotent mathematics is evident from the wide range of problems it can address. For a comprehensive overview of this  field and its applications, readers may refer to the survey article by Litvinov \cite{Litv} and the extensive bibliography provided therein.

Such fundamental  mathematical notion as convexity also found its idempotent analogues.  Max-plus convex sets were introduced in \cite{Z}.
Max-min convexity was studied in \cite{NS} and \cite{NS1}. The B-convexity based on the operations of the maximum and the multiplication was studied in \cite{BCh}.

Let us remark that many known abstract convexity structures have categorical nature and are  generated by monads \cite{R1}.
Topological and categorical properties of the functor of max-plus idempotent probabilistic measures  were studied in \cite{Zar}. In particular, the idempotent measure monad was constructed. The monad of $\cdot$-measures (functionals which preserve the maximum and the multiplication  operations) was introduced in  \cite{R6}.  The convexity generated by this monads coincide with the convexities based on the maximum and the multiplication  operations and were described in \cite{R2} in more general context.

Although idempotent measures are not additive and corresponding functionals are not linear, there are some parallels between topological properties of the functor of probability measures and the functor of idempotent probability measures (see for example \cite{Zar}, \cite{BR} and \cite{Radul}) which are based on existence of  natural equiconnectedness structure on both functors. However, some differences appeared when topological properties of the barycenter map were studying, because of differences of linear convexity and max-plus convexity.
The problem of the openness of the barycenter map of probability measures was investigated  in \cite{Fed}, \cite{Fed1},  \cite{Eif}, \cite{OBr} and \cite{Pap}. In particular, it is proved in \cite{OBr} that the barycenter map for a compact convex set  in a locally convex space is open if and only if the map $(x, y)\mapsto 1/2 (x + y)$ is open.
Zarichnyi defined in \cite{Zar} the idempotent barycenter map for max-plus idempotent probability measures and asked if the analogous characterization is true. A negative answer to this question was given  in  \cite{Radul1} where also another characterization of openness of idempotent barycenter map was given.

The problem of the softness of the barycenter map of probability measures was investigated  in \cite{Fed1}, \cite{Radul2} and \cite{Radul3}. In particular, it was proved in \cite{Fed1} that the product of $\omega_1$ barycentrically open metrizable convex compacta is barycentrically soft. Instead, it was proved in \cite{Radul4} that for the Tychonov cube of weight $\omega_1$ the idempotent barycenter map is not  soft. The problem when idempotent barycenter map  is a trivial bundle with the fiber Hilbert cube was studied in \cite{Radul5}.

W.Briec asked if it is possible to obtain some analogues of mentioned   results (concerning max-plus semiring) for another idempotent models. The idempotent barycenter map for $\ast$-measures was introduced in \cite{R2} where $\ast$ is any continuous t-norm. We investigate topological properties of that map when t-norm is the multiplication operation. To achieve these goals, we prove that the max-plus convexity is isomorphic to the convexity based on the maximum and multiplication operations.

\section{Idempotent measures and convexities} In what follows, all spaces are assumed to be compacta (compact Hausdorff space) except for $\R$ and the spaces of continuous functions on a compactum. All maps are assumed to be continuous.  We  denote by $C(X)$ the
Banach space of continuous functions on a compactum  $X$ endowed with the sup-norm.  We also consider the natural lattice operations $\vee$ and $\wedge$ on $C(X)$ and  its sublattices $C(X,[0,+\infty))$ and $C(X,[0,1])$. For any $c\in\R$ we shall denote the
constant function on $X$ taking the value $c$ by $c_X$.

\begin{df}\cite{Zar} A functional $\mu: C(X) \to \R$ is called an idempotent   measure (a Maslov measure) if

\begin{enumerate}
\item $\mu(1_X)=1$;
\item $\mu(\lambda_X+\varphi)=\lambda+\mu(\varphi)$ for each $\lambda\in\R$ and $\varphi\in C(X)$;
\item $\mu(\psi\vee\varphi)=\mu(\psi)\vee\mu(\varphi)$ for each $\psi$, $\varphi\in C(X)$.
\end{enumerate}

\end{df}

Let us remark that that the term 'idempotent  measure' is referred to as the 'idempotent probability measure' in \cite{Zar}. However, we adopt a simplified terminology to prevent confusion with the concept of a probability measure. In our framework, we equate a probability measure with a linear, positively defined normed functional, as established by Riesz Theorem. Consequently, the concept of an idempotent measure serves as an idempotent analogue of a probability measure. This analogy arises by substituting the usual addition and multiplication operations with the maximum and addition operations, respectively.

Let $IX$ denote the set of all idempotent  measures on a compactum $X$. We consider topologically
$IX$ as a subspace of $\R^{C(X)}$.



Let $\R_{\max}=\R\cup\{-\infty\}$ be the metric space endowed with the metric $\varrho$ defined by $\varrho(x, y) = |e^x-e^y|$.
 The convention $-\infty + x=-\infty$ and $-\infty \vee x=x$ allows us to extend the operations $\vee$ and $+$  over $\R_{\max}$. By $[-\infty,0]$ we denote the corresponding subset of $\R_{\max}$.

Max-plus convex sets were introduced in \cite{Z} and found many applications. Let $\tau$ be a cardinal number. Given $x, y \in \R^\tau$ and $\lambda\in[-\infty,0]$, we denote by $y\vee x$ the coordinatewise
maximum of x and y and by $\lambda+ x$ the vector obtained from $x$ by adding $\lambda$ to each of its coordinates. A subset $A$ in $\R^\tau$ is said to be  max-plus convex if $(\alpha+ a)\oplus  b\in A$ for all $a, b\in A$ and $\alpha\in[-\infty,0]$. It is easy to check that $A$  is   max-plus convex iff $\vee_{i=1}^n\lambda_i+ x_i\in A$ for all $x_1,\dots, x_n\in A$ and $\lambda_1,\dots,\lambda_n\in[-\infty,0]$ such that $\vee_{i=1}^n\lambda_i=0$. In the following by max-plus convex compactum we mean a max-plus convex compact subset of $\R^\tau$. It is shown in \cite{Zar} that $IX$ is a compact max-plus subset of $\R^{C(X)}$.

\begin{df}\cite{R6} A functional $\mu: C(X,[0,1]) \to [0,1]$ is called an $\cdot$ -measure  if

\begin{enumerate}
\item $\mu(1_X)=1$;
\item $\mu(\lambda_X\cdot\varphi)=\lambda\cdot\mu(\varphi)$ for each $\lambda\in[0,1]$ and $\varphi\in C(X,[0,1])$;
\item $\mu(\psi\vee\varphi)=\mu(\psi)\vee\mu(\varphi)$ for each $\psi$, $\varphi\in C(X)$.
\end{enumerate}

\end{df}

Let us denote that the term $\cdot$ -measure is not used in \cite{R6}. It is a particular case of the term $\ast$-measure from \cite{Sukh} where $\ast$ denote any continuous t-norm.

We denote by $A^\cdot(X)$) the space of all $\cdot$ -measures considered as a subspace of  the space $[0,1]^{C(X,[0,1])}$ with the product topology. 

 By $\delta_{x}$ we denote the Dirac measure supported by the point $x\in X$. Let $K$ be a compactum. Consider the subset $A^\cdot_\omega(K)\subset A^\cdot(K)$ defined as follows $$A^\cdot_\omega(K)=\{\bigvee_{i=1}^n\lambda_i\cdot \delta_{x_i}\mid n\in\N,\text{ }x_1,\dots, x_n\in A \text{ and }\lambda_1,\dots,\lambda_n\in[0,1]$$ $$\text{ such that }\bigvee_{i=1}^n\lambda_i=1\}.$$
It is known that $A^\cdot_\omega(K)$ is dense in $A^\cdot(K)$ \cite{Sukh}.

Let $T$ be a set. Given $x, y \in [0,1]^T$ and $\lambda\in[0,1]$, we denote by $y\vee x$ the coordinatewise
maximum of $x$ and $y$ and by $\lambda\ast x$ the point with coordinates $\lambda\ast x_t$, $t\in T$. We follow the convention that the operation $\cdot$ has higher priority than the maximum. A subset $A$ in $[0,1]^T$ is said to be  max- $\cdot$ convex if $\lambda\cdot a\vee  b\in A$ for all $a, b\in A$ and $\lambda\in[0,1]$. It is easy to check that $A$  is   max- $\cdot$ convex iff $\bigvee_{i=1}^n\lambda_i\cdot x_i\in A$ for all $x_1,\dots, x_n\in A$ and $\lambda_1,\dots,\lambda_n\in[0,1]$ such that $\bigvee_{i=1}^n\lambda_i=1$. In the following by max- $\cdot$ convex compactum we mean a max- $\cdot$ convex compact subset of $[0,1]^T$. The fact  that the set $ A^\cdot(K)$ is a max- $\cdot$ convex compact subset of $[0,1]^{C(X,[0,1])}$ is a particular case of a result from \cite{Sukh} obtained for any continuous t-norm.

The main goal of this section is to establish some correspondence between max-plus convex and   max- $\cdot$ convex compacta which preserves respective operations. Let us remark that it would be technically easier if we considered max-plus convex compacta as subsets of $\R_{\max}^\tau$. But many results are obtained for max-plus convex compacta as subsets of $\R^\tau$. Moreover, since idempotent measures are defined on functions from $C(X)$, we have a problem to define barycenter map for  max-plus convex compacta as subsets of $\R_{\max}^\tau$ without changing definition of the idempotent measure. So, we will do necessary technical work in this paper which in particular demonstrate that both approaches are equivalent.

It is also worth to notice that an isomorphism between max-plus convexity in $[-\infty,+\infty)^n$ and max- $\cdot$ convexity in $[0,+\infty)^n$ for a natural $n$ was mentioned without details in \cite{BCh}.

Let $S$ be a max- $\cdot$ convex compactum. A continuous map $h:S\to\R^\tau$ is called  $(\cdot,+)$-affine if for each $\lambda\in [0,1]$, $s$, $k\in S$ we have $h(\lambda\cdot s\vee  k)=(\ln(\lambda)+h(s))\vee  h(k)$ (we put $\ln(0)=-\infty$ and $\exp(-\infty)=0$). It is easy to check that $h(S)$ is a max-plus convex compactum.

Conversely, let $K$ be a max-plus convex compactum. A continuous map $g:K\to[0,1]^\tau$ is called an $(+,\cdot)$-affine if for each $\lambda\in [-\infty,0]$, $s$, $k\in K$ we have $g((\lambda+ s)\vee  k)=\exp(\lambda)\cdot g(s)\vee  g(k)$.  It is easy to check that $g(K)$ is a max- $\cdot$ convex compactum.

If a $(\cdot,+)$-affine map $h:S\to h(S)$ is a homeomorphism, we call it $(\cdot,+)$-isomorphism. It is easy to check that the inverse  $h^{-1}:h(S)\to S$ is $(+,\cdot)$-affine and we call it $(+,\cdot)$-isomorphism.

For $n\in\N$ let us define the map $f_n:[-\infty,0]\to \R$ as follows $$f_n(t)=\begin{cases}
t,&t>-n,\\
-n,&t\le-n\end{cases}
$$

for $t\in[-\infty,0]$. Evidently, $f_n$ is continuous.

\begin{proposition}\label{af} We have $f_n((\lambda+ s)\vee  k)=(\lambda+f_n(s))\vee  f_n(k)$ for each $\lambda\in [-\infty,0]$, $s$, $k\in [-\infty,0]$.
\end{proposition}


Consider the diagonal product $f=(f_i)_{i\in\N}:[-\infty,0]\to \R^\N$ of the maps $f_n$.

\begin{proposition}\label{emaf} The map $f$ is an embedding such that $f((\lambda+ s)\vee  k)=(\lambda+f(s))\vee  f(k)$ for each $\lambda\in [-\infty,0]$, $s$, $k\in [-\infty,0]$.
\end{proposition}


Composing the map $f$ with the $\ln$, we obtain the map $g:[0,1]\to \R^\N$ which is a $(\cdot,+)$-affine embedding.
Let $T$ be a set. Taking product of $T$ copies of the map $g$ we obtain the map $g_T:[0,1]^T\to \R^{\N\times T}$ which is a $(\cdot,+)$-affine embedding. Restricting the map  $g_T$ to any max- $\cdot$ convex compactum  we obtain the following theorem.

\begin{theorem}\label{genemaf} Let $S$ be a max- $\cdot$ convex compactum. There exists a $(\cdot,+)$-affine embedding $h:S\to\R^\tau$.
\end{theorem}

\section{Idempotent barycenter mappings}

Let $A\subset  \R^T$ be a compact max-plus convex subset. For each $t\in T$ we put $f_t=\pr_t|_A:A\to \R$ where $\pr_t:\R^T\to\R$ is the natural projection.    Given $\mu\in IA$, the point $\beta_A(\mu)\in\R^T$ is defined by the conditions $\pr_t(\beta_A(\mu))=\mu(f_t)$ for each $t\in T$. It is shown in \cite{Zar} that $\beta_A(\mu)\in  A$ for each $\mu\in I(A)$ and the map $\beta_A : I(A)\to A$ is continuous.
The map $\beta_A$ is called the idempotent barycenter map. It follows from results of \cite{Zar} that for each compactum $X$ we have $\beta_{IX}\circ I(\delta X)=\id_{IX}$ and for each map $f:X\to Y$ between compacta $X$ and $Y$ we have $\beta_{IY}\circ I^2f=If\circ\beta_{IX}$.

Let $K\subset  [0,1]^T$ be a compact max- $\cdot$ convex subset. For each $t\in T$ we put $f_t=\pr_t|_K:K\to [0,1]$ where $\pr_t:[0,1]^T\to[0,1]$ is the natural projection.    Given $\mu\in A^\cdot(K)$, the point $\beta^\cdot_K(\mu)\in[0,1]^T$ is defined by the conditions $\pr_t(\beta^\cdot_K(\mu))=\mu(f_t)$ for each $t\in T$. It follows from results of \cite{R2} that $\beta^\cdot_K(\mu)\in  K$ for each $\mu\in A^\cdot(K)$ and the map $\beta^\cdot_K : A^\cdot(K)\to K$ is continuous. We call this map max- $\cdot$ barycenter map.

We call two continuous maps $f:X\to Y$ and $g:T\to Z$ homeomorphic if there exists a pair of homeomorphisms $l:X\to T$ and $e:Y\to Z$ such that $e\circ f=g\circ l$. The main goal of this section is to show  that  the maps $\beta^\cdot_K : A^\cdot(K)\to K$ and $\beta_h(K) : I(h(K))\to h(K)$ are homeomorphic where  $h:K\to\R^\tau$ is the $(\cdot,+)$-affine embedding from Theorem \ref{genemaf}. This fact allow us to use results from \cite{Radul1}, \cite{Radul4} and \cite{Radul5} about topological properties of idempotent barycenter map for studying properties of max- $\cdot$ barycenter map.

Let $K\subset  [0,1]^T$ be a compact max- $\cdot$ compactum and $h:K\to\R^\tau$ is a $(\cdot,+)$-affine embedding. Denote $S=h(K)$.  Define the map $l_h:A^\cdot(K)\to IS$ as follows $$l_h(\nu)(\psi)=\max\{\psi(h(x))+\inf\{\ln(\nu(\varphi))\mid \varphi\in C(K,[0,1] \text{ such that } \varphi(x)=0\}\mid x\in X\}$$ for $\nu\in A^\cdot(K)$ and $\psi\in C(S)$.

\begin{lemma}\label{l} We have $l_h(\bigvee_{i=1}^n\lambda_i\cdot \delta_{x_i})=\bigvee_{i=1}^n(\ln(\lambda_i)+ \delta_{h(x_i)})$ for each $\bigvee_{i=1}^n\lambda_i\cdot \delta_{x_i}\in A^\cdot_\omega(K)$.
\end{lemma}


\begin{theorem}\label{hombar} The pair $(l_h,h)$ is a homeomorphism of maps $\beta^\cdot_K : A^\cdot(K)\to K$ and $\beta_S : I(S)\to S$.
\end{theorem}


\section{Topological properties of the max- $\cdot$ barycenter map}

We start with investigation of openness of the max- $\cdot$ barycenter map using the homeomorphism described in the previous section and the results from \cite{Radul1} about openness of the the idempotent barycenter map. It is easy to check that any homeomorphism of maps preserves openness.

\begin{theorem}\label{freeopen} The map $\beta^\cdot_{A^\cdot X} : A^\cdot(A^\cdot X)\to A^\cdot X$ is open for each compactum $X$.
\end{theorem}

\begin{proof} It was shown in \cite{Radul1} that the  map $\beta_{I X} : I(I X)\to I X$ is open for each compactum $X$. Hence sufficient demonstrate that the maps $\beta^\cdot_{A^\cdot X} : A^\cdot(A^\cdot X)\to A^\cdot X$ and $\beta_{I X} : I(I X)\to I X$ are homeomorphic. To achieve that, it is enough to find  a $(\cdot,+)$-affine homeomorphism $g_X:A^\cdot X\to I X$. We put $$g_X(\nu)(\psi)=\max\{\psi(x)+\inf\{\ln(\nu(\varphi))\mid \varphi\in C(X) \text{ such that } \varphi(x)=0\}\mid x\in X\}$$ for $\nu\in A^\cdot X$ and $\psi\in C(X)$.
\end{proof}

Put $J=\{(t,p)\in [0,1]\times [0,1]\mid t\vee p=1\}$. Let $X$ be a  max- $\cdot$ convex compactum.  We consider the map of max- $\cdot$ convex combination $s_X:X\times X\times J\to X$ defined by the formula $s_X(x,y,t,p)=t\cdot x\vee p\cdot y$. The map of max-plus convex combination $p_Y:Y\times Y\times J_0\to Y$, where $J_0=\{(t,p)\in [-\infty,0]\times [-\infty,0]\mid t\oplus p=0\}$ and  $Y$ is a  max-plus convex compactum, was introduced in \cite{Radul1}. Define the map $l:J\to J_0$ as follows $l(t,p)=(\ln t,\ln p)$. It is easy to see that for any $(\cdot,+)$-affine homeomorphism $h:X\to Y$ the pair $(h\times h\times l,h)$ is a homeomorphism of maps $s_X$ and $p_Y$. Using previous results and  this homeomorphism we obtain the max- $\cdot$ analogue of  Theorem 4.4 from \cite{Radul1}.

\begin{theorem}\label{main} Let $X$ be a max- $\cdot$ convex compactum. Then the following statements are equivalent:
\begin{enumerate}
\item the map $\beta^\cdot_X|_{A^\cdot_\omega X}:A^\cdot_\omega X\to X$ is open;
\item the map $\beta^\cdot_X$ is open;
\item the map $s_X$ is open.
\end{enumerate}
\end{theorem}

By $D=\{0,1\}$ we denote a two-point discrete compactum and by $K$ the max-plus convex compactum $ID$. It was shown in \cite{Radul1} that the map $\beta_K$ is open but the map $(x, y)\mapsto x \vee y$ is not. Since the maximum operation $\vee$ is an idempotent analogue of the addition operation, this result gives a negative answer to the question stated by Zarichnyi  \cite{Zar} and demonstrate some difference  between the theory of probability measures and idempotent measures. Let us remark that the map $\beta^\cdot_{A^\cdot D}:A^\cdot(A^\cdot D)\to A^\cdot D$ is also open by Theorem \ref{freeopen}, what demonstrates that we can not substantially simplify the characterization given in the previous theorem.

It is obvious that a product of a family of  (max-plus, max- $\cdot$)  convex compacta has a natural structure of (max-plus, max- $\cdot$) convexity with coordinatewise  operations.
It is proved in \cite{Fed1} that the product of barycentrically open compact convex sets (i.e. compact convex
sets for which the barycentre map  is open) is again barycentrically open. The situation is different for the idempotent barycentre map. It is shown in \cite{Radul1} that the idempotent barycenter map $\beta_S:IS\to S$ is not open for the max-plus  convex compactum $S=ID\times ID$, although the map $\beta_{I D}$ is open. We can translate this statement to the max- $\cdot$ model by means of the  $(\cdot,+)$-affine homeomorphism $g_D\times g_D:A^\cdot D\times A^\cdot D\to ID\times ID$, where $g_D$ is the homeomorphism introduced in the proof of Theorem \ref{freeopen}.

\begin{proposition} The  map $\beta^\cdot_K:A^\cdot K\to K$ is not open for the max- $\cdot$  convex compactum $K=A^\cdot D\times A^\cdot D$.
\end{proposition}

Let us remark that the map $\beta^\cdot_{A^\cdot D}$ is  open by Theorem \ref{freeopen}.

However, it is shown in in \cite{Radul1} that the idempotent barycenter map $\beta_P:IP\to P$ is  open where $P$ is the product of intervals $P=[a,b]^\tau$, $a$, $b\in\R$. The analogous result we have for  max- $\cdot$ barycenter map.

\begin{proposition}\label{prodopen} Let $a$, $b\in(0,1]$ and $P=[a,b]^\tau$. The  map $\beta^\cdot_P:A^\cdot P\to P$ is  open.
\end{proposition}

\begin{proof} Consider the map $\ln:[a,b]\to\R$ and its power $h=(\ln)^\tau:[a,b]^\tau\to\R^\tau$. It is easy to see that the map $h$ is a $(\cdot,+)$-affine embedding. By Theorem \ref{hombar} the pair $(l_h,h)$ is a homeomorphism of the maps $\beta^\cdot_P:A^\cdot P\to P$ and $\beta_K:IK\to K$ where $K=[\ln a,\ln b]^\tau$ is an max-plus  convex subset in $\R^\tau$. Thus, openness of the map $\beta_K:IK\to K$ implies openness of the map $\beta^\cdot_P:A^\cdot P\to P$.
\end{proof}

The next proposition demonstrates that we can not use the previous  arguments for the cube $[0,1]^\tau$.

\begin{proposition} There is no $(\cdot,+)$-affine embedding of $[0,1]$ into $\R$.
\end{proposition}

\begin{proof} Suppose the contrary. Then there exists a $(\cdot,+)$-affine embedding of $s:[0,1]\to \R$. Put $s(0)=a$ and $s(1)=b$. Since $s$ is monotone, we have $a<b$. Choose $\lambda>0$ such that $\ln\lambda<a-b$. Then we have $a<s(\lambda)=s(\lambda\cdot 1\vee0)=(\ln\lambda+b)\vee a\le a\vee a=a$ and we obtain a contradiction.
\end{proof}

So, we have the following problem.

\begin{question} Is the map $\beta^\cdot_{[0,1]^\tau}:A^\cdot ([0,1]^\tau)\to [0,1]^\tau$ open?
\end{question}

Now, let's turn our attention to a stricter class of maps -- specifically, soft maps.  A map $f:X\to Y$ is said to be (0-)soft if for any (0-dimensional) paracompact space $Z$, any closed subspace $A$ of $Z$ and
maps $\Phi:A\to X$ and $\Psi:Z\to Y$ with $\Psi|A=f\circ\Phi$ there exists a
map $G:Z\to X$ such that $G|A=\Phi$ and $\Psi=f\circ G$. This notion is
introduced by E.V.Shchepin \cite{Shchep1}. Let us remark that each 0-soft map is open and 0-softness is equivalent to the openness for all the maps between metrizable compacta.

It was proved in \cite{Radul4} that softness of the map $\beta_{IX}$ implies metrizability of the compactum $X$. It is easy to check that softness of maps is  preserved by  homeomorphisms between maps. Using the $(\cdot,+)$-affine homeomorphism $g_X:A^\cdot X\to I X$ introduced in the proof of Theorem \ref{freeopen}, we obtain the following statement.

\begin{theorem} Let $X$ be a  compactum such that the map $\beta^\cdot_{A^\cdot X} : A^\cdot(A^\cdot X)\to A^\cdot X$ is soft. Then $X$ is metrizable.
\end{theorem}

Let us remark that we have the analogous result for the classic barycenter map of probability measures \cite{Radul2}. The difference between classic and idempotent cases appears when we discuss softness of barycenter maps for products. Fedorchuk proved in \cite{Fed1} that  each product of $\omega_1$ barycentrically open convex metrizable compacta  is barycentrically soft. It is proved  in \cite{Radul4} that the idempotent barycenter map  $\beta_{[0,1]^{\omega_1}}$ is not soft. It is easy to modify this result in order to obtain that the idempotent barycenter map  $\beta_{[-1,0]^{\omega_1}}$ is not soft. Applying arguments as in the proof of \ref{prodopen}, we obtain that the max- $\cdot$ barycenter map  $\beta^\cdot_{[\exp(-1),1]^{\omega_1}}$ is not soft. The problem of existing of max-plus convex compactum with soft idempotent barycenter map is posed in \cite{Radul4}. The same problem is open for the max- $\cdot$ barycenter map.


\begin{thebibliography}{}



\bibitem{BR} T.Banakh, T.Radul, {\em F-Dugundji spaces, F-Milutin spaces and absolute F-valued retracts}, Topology Appl. {\bf 179} (2015), 34--50.

\bibitem{BCh}  W.Briec, Ch.Horvath {\em Nash points, Ky Fan inequality and equilibria of abstract economies in Max-Plus and $\mathbb B$-convexity,} J. Math. Anal. Appl. {\bf 341} (2008), 188--199.



\bibitem{Eif} L.Q. Eifler, {\em Openness of convex averaging}, Glasnik Mat. Ser. III, {\bf 32} (1977), no. 1, 67--72.


\bibitem{Fed} V.V.~Fedorchuk, {\em  On a barycentric map of probability measures}, Vestn. Mosk. Univ, Ser. I, No 1, (1992), 42--47.

\bibitem{Fed1} V.V.~Fedorchuk, {\em On barycentrically open bicompacta}, Siberian Mathematical Journal, {\bf 33} (1992), 1135--1139.


\bibitem{Litv} G. L.Litvinov, {\em The Maslov dequantization, idempotent and tropical mathematics: a very brief
introduction}, Idempotent mathematics and mathematical physics, 1--17, Contemp. Math., 377, Amer. Math. Soc., Providence, RI, 2005.

\bibitem{NS} V. Nitica, I. Singer, {\em Contributions to max-min convex geometry. I. Segments}, Linear Algebra Appl. {\bf 427}, N7 (2008) 1439--1459.

\bibitem{NS1} V. Nitica, I. Singer, {\em Contributions to max-min convex geometry. II. Semispaces
and convex sets}, Linear Algebra Appl. {\bf 428}, N8-9 (2008) 2085--2115.

\bibitem{OBr} R.C. O'Brien, {\em On the openness of the barycentre map}, Math. Ann., {\bf 223} (1976), 207--212.

\bibitem{Pap} S. Papadopoulou, {\em  On the geometry of stable compact convex sets}, Math. Ann., {\bf 229} (1977), 193--200.

\bibitem{R1} T.Radul, {\em Convexities generated by L-monads},  Applied Categorical Structures {\bf 19} (2011) 729--739.

\bibitem{R2} T.Radul, {\em On t-normed integrals with respect to possibility capacities on compacta}, Preprint,  arXiv:2111.06612.

\bibitem{Radul} T.~Radul, {\em Absolute retracts and equiconnected monads}, Topology Appl. {\bf 202} (2016), 1--6.

\bibitem{Radul1} T.~Radul, {\em On the openness of the idempotent barycenter map}, Topology Appl.{\bf 265} (2019) 106809.

\bibitem{Radul2} T.~Radul, {\em On the baricentric map of probability measures}, Vestn. Mosk. Univ., Ser. I  (1994), No.1, 3--6.

\bibitem{Radul3} T.~Radul, {\em On baricentrically soft compacta}, Fund.Math. {\bf 148} (1995), 27--33.

\bibitem{Radul4} T.~Radul, {\em Idempotent measures:absolute retracts and soft maps}, Topological Methods in Nonlinear Analysis {\bf 59}, N 2B, (2022), 1029--1045.

\bibitem{Radul5} T.~Radul, {\em Fibration of idempotent measures}, Ukr.Math.Journal {\bf 72}, N 11, (2021), 1784--1793.


\bibitem{R6} T.Radul, {\em A functional representation of the capacity multiplication monad},  Visnyk of the Lviv Univ. Series Mech. Math. {\bf 86} (2018) 125--133.

\bibitem{Sukh} Kh. Sukhorukova, {\em Spaces of non-additive measures generated by triangular norms},   Mat. Stud. {\bf 59} (2023) 215--224.

\bibitem{Shchep1} E.V.Shchepin, {\em Topology of limit spaces of uncountable inverse spectra}, Russian Mathematical Surveys {\bf 31} (1976), 155--191.



\bibitem{Zar} M.~Zarichnyi, {\em Spaces and mappings of idempotent measures}, Izv. Ross. Akad. Nauk Ser. Mat. {\bf 74} (2010), 45--64.

\bibitem{Z} K. Zimmermann, {\em A general separation theorem in extremal algebras}, Ekon.-Mat. Obz. {\bf 13} (1977) 179--201.

\end{thebibliography}
\end{document}